\documentclass[a4paper,12pt]{paper}

\usepackage[T1]{fontenc}
\usepackage[utf8]{inputenc}
\usepackage[english]{babel}
\usepackage{amsthm,amsmath,amsfonts,amssymb}

\usepackage[scaled=.92]{helvet}

\usepackage[]{sfmath}

\usepackage{tikz}
\usetikzlibrary{shapes,arrows}
\usepackage{url}
\usepackage{listings}
\usepackage{algorithm2e}
\usepackage{algorithmic}

\usepackage[abs]{overpic}
\numberwithin{figure}{section}
\numberwithin{equation}{section}

\newcommand{\R}{\mathbb{R}}

\def\bei#1{\vrule width 0.4pt height 14pt depth 9pt
           \lower 8pt \hbox{$ _{\hbox{}\, #1}$}\!}
\renewcommand{\div}[1][]{{\operatorname{div}_{#1} \,}} %divergence operator with extra stuff for tangential divergence
\renewcommand{\d}[0]{{\operatorname{d}}}

\newcommand{\velo}{u}
\newcommand{\p}{p}

\newcommand{\dvis}{\mu}
\newcommand{\dens}{\rho}
\newcommand{\n}{n}
\newcommand{\V}{V}

\begin{document}

\title{A Two Stage CVT / Eikonal Convection Mesh Deformation Approach for Large Nodal Deformations}

\author{Stephan Schmidt\thanks{Dr. Stephan Schmidt, University of Würzburg, Germany; stephan.schmidt@mathematik.uni-wuerzburg.de}}

\maketitle

\begin{abstract}
A two step mesh deformation approach for large nodal deformations, typically arising from non-parametric shape optimization, fluid-structure interaction or computer graphics, is considered. Two major difficulties, collapsed cells and an undesirable parameterization, are overcome by considering a special form of ray tracing paired with a centroid Voronoi reparameterization. The ray direction is computed by solving an Eikonal equation. With respect to the Hadamard form of the shape derivative, both steps are within the kernel of the objective and have no negative impact on the minimizer. The paper concludes with applications in 2D and 3D fluid dynamics and automatic code generation and manages to solve these problems without any remeshing. The methodology is available as a FEniCS shape optimization add-on at~\url{http://www.mathematik.uni-wuerzburg.de/~schmidt/femorph}.
\end{abstract}

\section{Introduction}
Shape optimization is a special sub-class of PDE constrained optimization. Contrary to more common problems of minimizing a functional, shape optimization naturally leads to the question of how to parameterize the unknown domain and the respective derivative or shape update. Typically, these two problems are treated separately by first defining a parameterization, usually only of the deformation but not the shape itself, and then computing the derivative with respect to the parameterization on a separate and independent mesh. Within aerospace applications, the so-called Hicks-Henne ansatz functions~\cite{hickshenne} are often used. Other examples are free-form deformation or CAD splines with their respective parameters as design unknowns.

The advantage of this approach is that the discretization of the state and adjoint equation is essentially independent of the discretization of the deformation. Furthermore, the search space is limited to the image space of the deformation, usually resulting in a guaranteed minimum regularity but also restricting the maximum design space. The downside is, however, that the mapping of a deformation parameter to the PDE solution needs to be incorporated into the differentiation process of the objective, which can be costly to compute and requires the derivative of both the parameterization as well as the derivative of the PDE solution process with respect to the input node positions. Both aspects, when taken together, are often called ``mesh sensitivities''. Although they can be computed efficiently~\cite{GaWaMoWi07, NielsenPark2006}, one often needs information not available with most commercial CAD tools. Furthermore, the approach of separating the parameterization from the discretization of the PDE makes exploiting the structure of shape optimization problems difficult, as post-discretization, the problem can typically be seen as using a functional as design unknown. Other recent approaches to overcome the meshing problem are based on immersed boundaries and phase fields~\cite{Hecht1}.

The idea of shape calculus~\cite{DelfourZolesio, ZolesioSokolowski} is to directly differentiate the objective with respect to the input domain. Given sufficient regularity, typically Lipschitz-continuity of the boundary, the directional derivative stemming from this approach can be cast into the Hadamard-form, i.e.\ be transformed into a boundary integral expression with tangential kernel. After discretization, a descent in the objective can then be achieved by moving boundary nodes in normal direction, thereby combining deformation parameterization and PDE discretization into one. Furthermore, this approach does not require any consideration of the mesh deformation within the derivative chain, making the resulting procedure quite attractive and efficient from a numerical standpoint, as only a scalar boundary integral kernel involving primal and dual state needs to be evaluated. The downside of this approach, however, is that higher order schemes require additional considerations~\cite{SchulzRiemannian} or at least some other type of approximative Newton scheme or alternate space to find the gradient in~\cite{ProtazA2008}, sometimes called Sobolev-gradient descent. Because tangential nodal movement is within the kernel of the Hadamard-form, the parameterization of the optimal shape of this approach is often quite uneven and heavily influenced by the node placement on the original shape, leading to the problem of redistributing surface nodes to guarantee a satisfactory discretization of the final shape. The situation is quite similar to curvature flow, for which schemes with automatic node redistribution can for example be found in~\cite{GarckeEquiA2011}. The idea there is to use a modified weak form to compute the curvature flow, which automatically leads to an equidistant node spacing on curves. Well-parameterized multi-phase fluid interfaces are for example discussed in~\cite{GarckeInterface2014}.

Because shape optimization usually is not intrinsic to the boundary, that is traces of quantities defined in the volume are often required to calculate the boundary flow, and the methodology needs to extend to the three dimensional situation, we consider coupling the discretization of the continuous boundary flow with a discrete surface mesh quality measure stemming from a Voronoi partitioning usually used within mesh generation but not mesh deformation. After a well-parameterized surface is found, we use a ray-free ray tracing approach based on solving two Eikonal equations to deform the volume mesh robustly. The methodology is tested on optimizing an obstacle in an incompressible Navier--Stokes fluid, which requires very large nodal deformations. The new methodology works without re-meshing, which is highly advantageous for novel optimization methods like One-Shot~\cite{SIGS2011, GS09}.

\section{Nodal Shape Optimization and the Hadamard Theorem}
\subsection{Problem Definition}
Shape optimization problems are a subclass of PDE-constrained optimization, where parts of the geometry, that is the domain $\Omega \subset D \subset \R^d$, are the design unknowns to be found. The set $D$ is typically called the ``hold-all''. Within the context of nodal deformation, three typical cases can be identified:
\begin{align}
 \min\limits_{\Omega} J_1(\Omega) :=& \int\limits_{\Omega} f_1(x) \ \d x\label{eq:prob1} \\
 \min\limits_{\Omega} J_2(\Omega) :=& \int\limits_{\partial \Omega} f_2(s) \ \d s\label{eq:prob2} \\
 \min\limits_{\Omega} J_3(\Omega) :=& \int\limits_{\partial \Omega} \langle f_3(s), n(s)\rangle\label{eq:prob3} \ \d s,
\end{align}
where $\partial \Omega$ is the boundary of $\Omega$ with outward facing normal $n$. The above notation gives rise to the notion of a ``non-parametric'' approach, because $\Omega$ is formally used directly as the optimization variable, whereas a parametric approach would consider $\Omega(\xi)$ to be defined by a possibly smooth parameterization with parameter $\xi$, typically splines. Other alternatives are treating $\Omega(\xi)$ locally explicitly as the graph of a function with parameter $\xi$ or implicitly as a level set.

\subsection{Directional Derivatives and the Hadamard Theorem}
There are several approaches to define a perturbed domain $\Omega_t$. A detailed overview, as well as an analysis of the shape differentiability and the resulting gradient expressions can be found in~\cite{DelfourZolesio, ZolesioSokolowski}. A very basic summary will be given here. Within the context of ``perturbation of identity'', a sufficiently smooth vector field $V: D \rightarrow \R^d$ is considered, such that a perturbed domain can be defined by
\begin{align}
 \Omega_{\epsilon}[V] := \{x + \epsilon V(x) : x \in \Omega\}.
\end{align}
Then the derivative of~\eqref{eq:prob1} - \eqref{eq:prob3} in direction $V$ is given by the one-sided limit
\begin{align}
 dJ_i (\Omega)[V] := \lim\limits_{\epsilon\rightarrow 0^+}\frac{J_i(\Omega_\epsilon[V]) - J(\Omega)}{\epsilon}
\end{align}
and the gradient for the cases~\eqref{eq:prob1}-\eqref{eq:prob3} above is given by
\begin{align}
  dJ_1(\Omega)[V] := &\int\limits_{\partial \Omega} \langle V,n\rangle \ f_1(s) \ \d s\label{eq:lin1}\\
  dJ_2(\Omega)[V] := &\int\limits_{\partial \Omega} \langle V,n\rangle \left(\frac{\partial f_2(s)}{\partial n(s)} + \kappa(s) f_2(s)\right) \ \d s\label{eq:lin2} \\
  dJ_3(\Omega)[V] := &\int\limits_{\partial \Omega} \langle V,n\rangle \ \div f_3(s) \ \d s, \label{eq:lin3}
\end{align}
where $\kappa$ denotes the mean curvature and $\div$ refers to the standard divergence operator, interpreted as the trace inside the boundary integral. The derivation of the last identity, Equation~\eqref{eq:lin3}, can for example be found in~\cite{PhDSchmidtShape}. Although it is also possible to work with volume formulations of the respective linearizations above, we will focus on the respective boundary expressions here, which typically requires stronger smoothness assumptions but can result in very fast numerical algorithms only operating on the boundary. Thus, one is free from including the mesh deformation chain into the derivative computation, although this is possible~\cite{GaWaMoWi07}. Furthermore, any tangential perturbation for $V$ results in a zero directional derivative. Both properties will be exploited within the mesh reparameterization approach discussed here, such that the repair step is always inside the kernel of the objective.

One achieves descent in the objective $J_i$ by iterating over sets in accordance to i.e.\ a steepest descent algorithm by updating
\begin{align}\label{eq:Update}
 \partial \Omega_{k+1} := \{s + g_i(s) n(s) : s \in \partial \Omega_k\},
\end{align}
where $g_i$ denotes the respective scalar integration kernel from~\eqref{eq:lin1}-\eqref{eq:lin3} above. If $\Omega$ is a tessellated domain, i.e.\ a finite element or finite volume mesh, then the above formula provides an algorithm of how to move the surface vertices into the vertex normal direction to achieve descent. However, the resulting point distribution based on the above algorithm is typically rather uneven and irregular. Also, there is no obvious strategy of how to deform the volume mesh and standard approaches quickly lead to compressed or collapsed cells due to the typically very local movement of the boundary. As part of this work, a two stage centroidal voronoi surface reparameterization coupled to an Eikonal convection mesh deformation will be studied.

\section{From Remeshing to CVT / Eikonal Mesh Deformation}
%\subsection{From Remeshing to Mesh Deformation}
The two typical problems with mesh based nodal shape optimization are shown in Figure~\ref{fig:DestroyedMesh}:
\begin{figure}[h!]
\begin{center}
  \includegraphics[width=0.9\textwidth]{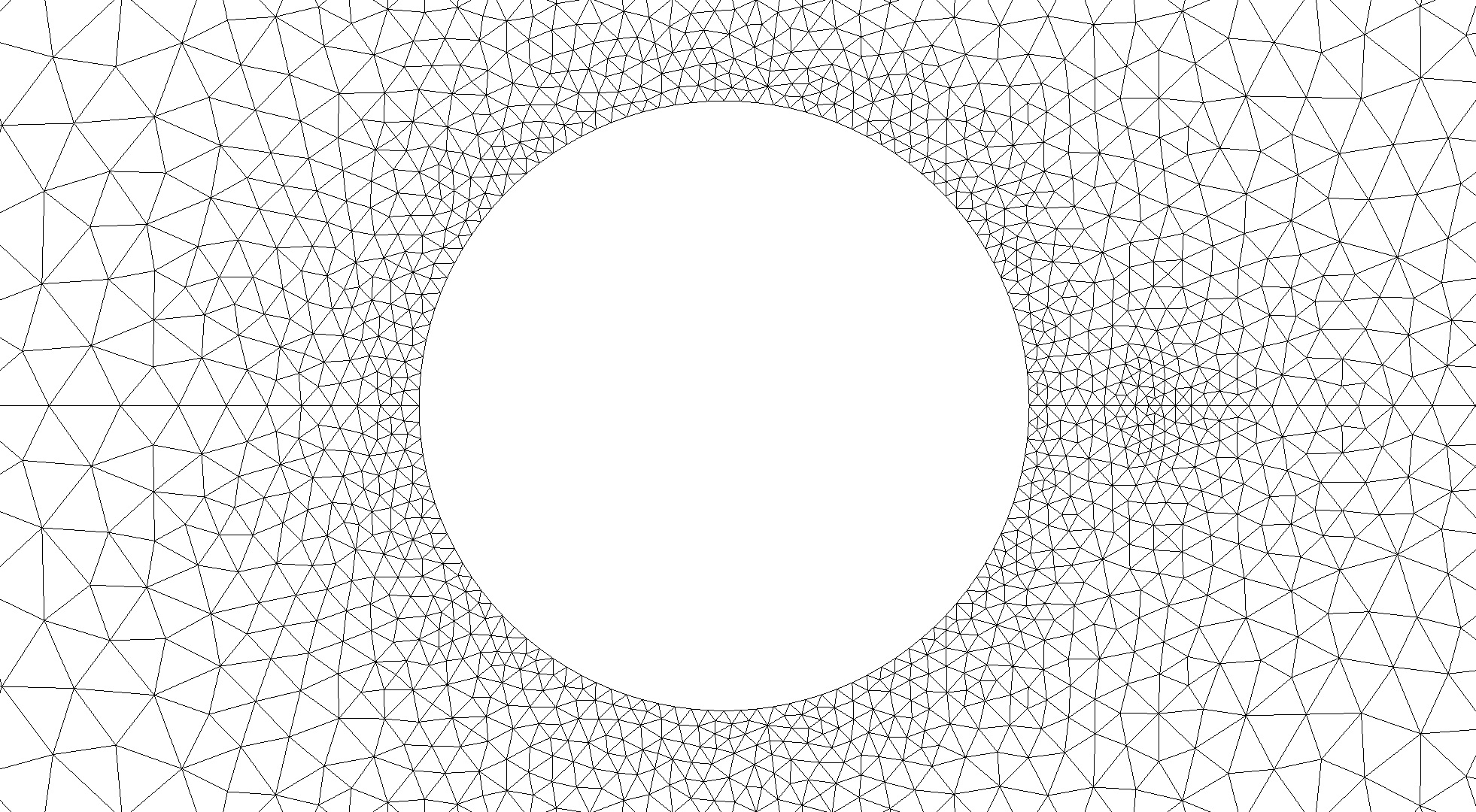}\vspace{0.1cm}\\
  \includegraphics[width=0.9\textwidth]{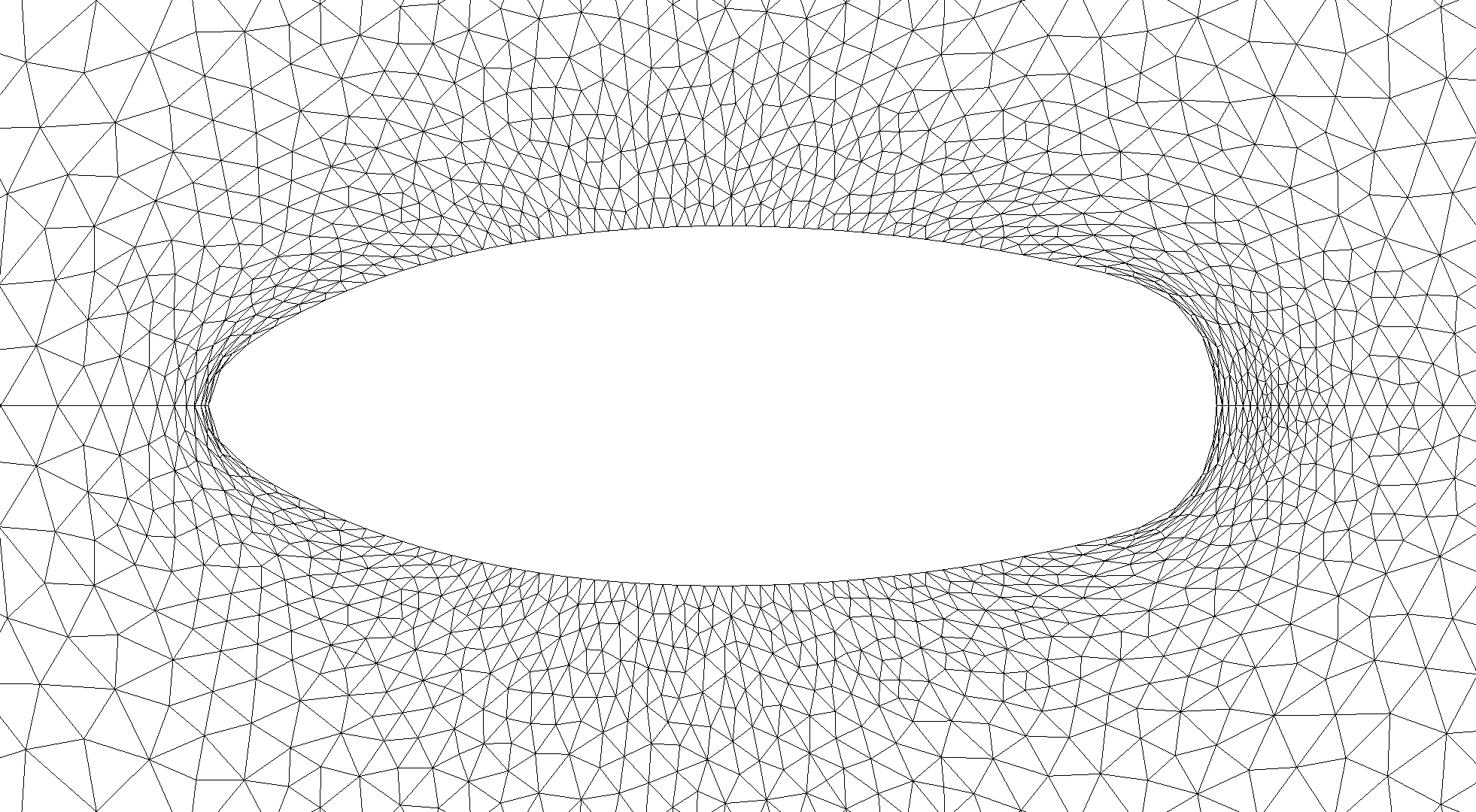}
\end{center}
  \caption{The initial mesh and the best optimum achievable when using a naive descent based on equation~\eqref{eq:Update}. Note the very coarse surface parameterization at the front and rear as well as the collapsed cells in the volume.}\label{fig:DestroyedMesh}
\end{figure}
Using an update based on equation~\eqref{eq:Update} naively results in a progressive deterioration of the parameterization during the optimization morph. Furthermore, a naive component-wise Laplacian mesh deformation does not prevent the occurrence of collapsed cells.  Each of these two problems will be addressed separately.

\subsection{Centroidal Voronoi Reparameterization}
\begin{figure}[h!]
\begin{center}
\usetikzlibrary{calc}
\begin{tikzpicture}[scale=3]
 \coordinate (Center) at (0.0, 0.0);
 \coordinate (x1) at (1.3, 0.0);
 \coordinate (x2) at (0.9, 0.8);
 \coordinate (v1) at (barycentric cs:Center=1.0,x1=1.0,x2=1.0);
 \coordinate (x3) at (-0.2, 0.9);
 \coordinate (v2) at (barycentric cs:Center=1.0,x2=1.0,x3=1.0);
 \coordinate (x4) at (-0.9, -0.1);
 \coordinate (v3) at (barycentric cs:Center=1.0,x3=1.0,x4=1.0);
 \coordinate (x5) at (0.2, -0.9);
 \coordinate (v4) at (barycentric cs:Center=1.0,x4=1.0,x5=1.0);
 \coordinate (v5) at (barycentric cs:Center=1.0,x5=1.0,x1=1.0);
 %\draw[help lines] (-4, -4) grid (4,4);
 %Draw primal mesh triangles
 \draw[thick] (Center) -- (x1) -- (x2) -- cycle;
 \draw[thick] (Center) -- (x2) -- (x3) -- cycle;
 \draw[thick] (Center) -- (x3) -- (x4) -- cycle;
 \draw[thick] (Center) -- (x4) -- (x5) -- cycle;
 %close
 \draw[thick] (Center) -- (x1) -- (x5) -- cycle;
 %Draw dual mesh
 \filldraw[thick,gray,dashed,opacity=0.5] (v1) -- (v2) --(v3) -- (v4) -- (v5) -- cycle;
 %Draw infinity edges
 \draw[thick,gray,dashed] (v1) -- ($(x1)!0.5!(x2)$);
 \draw[thick,gray,dashed] (v2) -- ($(x2)!0.5!(x3)$);
 \draw[thick,gray,dashed] (v3) -- ($(x3)!0.5!(x4)$);
 \draw[thick,gray,dashed] (v4) -- ($(x4)!0.5!(x5)$);
 \draw[thick,gray,dashed] (v5) -- ($(x5)!0.5!(x1)$);
 %Annotations
 \node at (0.1,-0.1) {$x_i$};
 \node at (0.5,0.25) {$W_i$};
 \node at (-0.55,0.67) {$\Omega_i$};
\end{tikzpicture}
\end{center}
\caption{Triangle patch $\Omega_i$ of the mesh in solid black lines. Dual cell $W_i$ based on triangle centroids around generator $x_i$ in gray. Primal mesh not Delaunay, because $x_i$ is not the centroid of the dual patch.}\label{fig:DualMesh}
\end{figure}
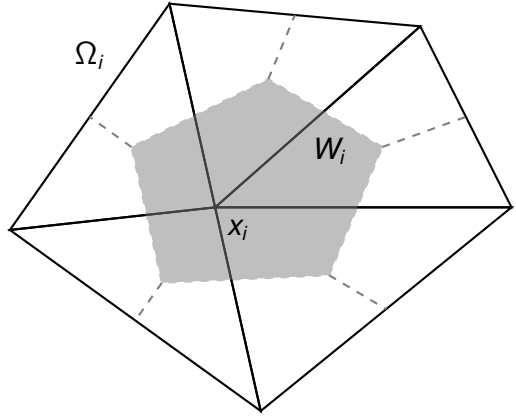

%According to~\cite{Chen2004_1} an optimal Delaunay triangulation $\mathcal{T}^*$ can be defined as the triangulation $\mathcal{T}$, for which a function $f \in C^{1}(\bar{\Omega})$ is best approximated by its piecewise linear and global continuous interpolation $f_\mathcal{T}$ based on $\mathcal{T}$, that is
%\begin{align*}
% Q(\mathcal{T}, f, p) := \vert f - f_\mathcal{T}\vert_{L^p(\Omega)}
%\end{align*}
%is minimized by $\mathcal{T}^*$ out of the set of all triangulations with at most $N$ vertices, the existence of which is given in~\cite{Chen2004_2} for convex $f$.
Since the main focus of this work is to achieve a robust mesh deformation for nodal shape optimization, we assume that an initial triangulated surface $\Gamma$ of the object to be optimized and a tessellation of the surrounding polygonal domain $\Omega$ is already given, usually by some form of Delaunay discretization. The dual of the Delaunay triangulation is the centroidal Voronoi tessellation. This gives rise to define the badness of a mesh consisting of $k$ vertices $\{x_i : i = 1,...,k\}$ by
\begin{align}\label{eq:MeshBadness}
 \mathcal{F}(\{x_i\}, W_i) := \sum_{i=1}^{k} \int\limits_{W_i} \rho(x) \Vert x - x_i \Vert^2_2 \ \d x,
\end{align}
where $\Omega = \cup_{i=1}^k W_i$ is a discretization of $\Omega$ into $k$ disjunct subsets. The set $\{x_i\}$ is also sometimes called the set of generators. Finally, $\rho$ is a varying density function to prescribe areas where a local clustering of points is desired. Typically, the sets $W_i$ are the dual cells of the triangulated primal mesh with vertices $x_i$, as shown in Figure~\ref{fig:DualMesh}.

It can be shown~\cite{CVT:ApplicationsAlgorithms} that a necessary condition for $\mathcal{F}$ to be minimized is that the dual patches $W_i$ are the Voronoi regions $W^*_i$
\begin{align*}
 W^*_i := \{x \in \Omega : \Vert x - x_i\Vert < \Vert x - x_j\Vert, j=1,...,k, i\neq j \}
\end{align*}
corresponding to $x_i$. Simultaneously, the vertices $x_i$ are the centroids of the corresponding $W_i$, because
\begin{align*}
 \nabla_{x_i} \mathcal{F}(\{x_i\}, W_i) &= 2 \int\limits_{W_i} \rho(x)\ \d x \cdot x_i - 2 \int\limits_{W_i} \rho(x)x \ \d x \stackrel{!}{=} 0 \\
 \Rightarrow x_i &= \frac{\int\limits_{W_i}\rho(x)x \ \d x}{\int\limits_{W_i} \rho(x) \ \d x},
\end{align*}
meaning choosing the centroid of $W_i$ for $x_i$ is a critical point if $W_i$ is considered fixed.
%\begin{align*}
% x_i = \frac{\int\limits_{W_i} \rho(x)x \d x}{\int\limits_{W_i} \rho(x) \ \d x}.
%\end{align*}
%When interpreting a triangular mesh as the dual of a not necessarily centroidal Voronoi tessellation, there is the correspondence that the Voronoi generators $x_i$ correspond to the mesh vertices, whereas the Voronoi regions $W_i$ correspond to the mesh edges, that is the connectivity.

Of equal importance to the node spacing is the node connectivity of the primal mesh. Vertices $x_i$ and $x_j$ of the primal mesh are connected by an edge if and only if the dual regions they generate are adjacent. Thus, dual patches $W_i$ in~\eqref{eq:MeshBadness} correspond to edges of the primal mesh. Because the goal here is to ensure a high quality mesh deformation during nodal shape optimization, a change in mesh connectivity is highly unwanted, corresponding to a remeshing procedure rather than a mesh deformation approach. Thus, we follow the notion of~\cite{Chen2004_1} by considering only the vertex positions, but not their connectivity, as the defining factor of mesh badness during nodal shape optimization. Contrary to the considerations before, patches $\Omega_i$ of the primal triangular mesh, rather than the Voronoi patches of the dual mesh, are used to define an alternate mesh badness measure. This is of considerable numerical advantage, because a costly computation of the Voronoi tessellation is omitted. The resulting mesh badness measure is thus given by
\begin{align}\label{eq:MeshBadness2}
 \mathcal{F}(\{x_i\}) := \sum_{i=1}^{k} \int\limits_{\Omega_i} \rho(x) \Vert x - x_i \Vert^2_2 \ \d x
\end{align}
and successively setting the mesh vertices to the critical points of~\eqref{eq:MeshBadness2} results in the CVT-I smoother as discussed in~\cite{Chen2004_1}.

\subsection{Surface Mesh Badness for Nodal Shape Optimization}
For the nodal shape optimization problem under consideration here, the above badness measure needs to be adapted, such that a vertex distribution that is a critical point of the badness measure is inside the kernel of the shape derivative. Consequently, the mesh repair step is a tangential movement of surface mesh nodes only and does not interfere with optimality, because mesh quality step and shape optimality/feasibility step are orthogonal to each other, due the shape derivatives~\eqref{eq:lin1}--\eqref{eq:lin3} being invariant under tangential deformation directions $V$.

A natural process when adapting~\eqref{eq:MeshBadness2} to a surface mesh would be to change the $\Vert.\Vert_2$-norm to the geodetic distance between points, as using the Euclidian distance will not result in the critical vertex set $\{x_i\}$ to be on the surface $\Gamma := \partial \Omega$. The computation of geodetic distances and their derivatives can be quite costly. However, the Euclidian distance within the local tangent plane to $\partial \Omega_i$ at $x_i$, that is $T_{x_i} \partial \Omega_i$, also results in critical points of~\eqref{eq:MeshBadness2} to stay within said tangent plane, as will be discussed next.

For any point $x$ in $\Omega_i\subset \R^d$, the projection to the tangent plane $T_{x_i} \Gamma$ is given by
\begin{align*}
 x_{\Gamma} := x - \langle x,n(x_i)\rangle \cdot n(x_i),
\end{align*}
where $n(x_i)$ is the normal at $x_i$. Using this definition, the surface mesh badness considered here is given by
\begin{align}\label{eq:MeshBadness3}
 \mathcal{F}(\{x_i\}) := \sum_{i=1}^{k_\Gamma} \int\limits_{\Gamma_i} \rho(s) \Vert (s - x_i)_{\Gamma} \Vert^2_2 \ \d s
\end{align}
and the critical set is given by points $x_i^*$ such that
\begin{align}\label{eq:Repair1}
\int\limits_{\Gamma_i} \rho(s)\left( -(s-x_i^*)_\Gamma - \langle s - x_i^*, n\rangle n'[x_i] (s-x_i^*)_\Gamma\right) \ \d s \stackrel{!}{=} 0,
\end{align}
where
\begin{align*}
 n'[x_i] := \nabla_{x} n(x) = \left(\frac{\partial}{\partial x}n(x_i)\right)^T \in \R^{d\times d}
\end{align*}
denotes the variation of the normal with respect to a shift of the origin point $x_i$ of the tangent plane. It can be shown~\cite{DelfourZolesio, PhDSchmidtShape, ZolesioSokolowski} that this object also lies within the tangent plane. The above relationship directly follows from
\begin{align*}
 &\nabla_{x_i} \Vert (s-x_i)_\Gamma\Vert^2_2\\
 =& \nabla_{x_i} \sum_{\ell=1}^d \left(s^\ell - x_i^\ell - \langle s - x_i, n\rangle n^\ell \right)^2\\
 =& \sum_{\ell=1}^d 2\left(s^\ell - x_i^\ell - \langle s - x_i, n\rangle n^\ell \right) \nabla_{x_i} \left[s^\ell - x_i^\ell - \langle s - x_i, n\rangle n^\ell \right] \\
 =& -2 (s-x_i)_\Gamma - 2\left\langle(s-x_i)_\Gamma, \nabla_{x_i} \left[\langle s-x_i,n\rangle\cdot n\right]\right\rangle,
\end{align*}
where upper indices denote vector components. Furthermore,
\begin{align*}
&\left\langle(s-x_i)_\Gamma, \nabla_{x_i} \left[\langle s-x_i,n\rangle\cdot n\right]\right\rangle \\
=&-\left\langle (s-x_i)_\Gamma, n\right\rangle \cdot n + \langle (s-x_i)_\Gamma, n\rangle n'[x_i] \cdot (s-x_i) + \left\langle s-x_i, n\right\rangle n'[x_i] \cdot (s-x_i)_\Gamma \\
=& \left\langle s-x_i, n\right\rangle n'[x_i] \cdot (s-x_i)_\Gamma,
\end{align*}
where the relationships $(a-b)_\Gamma = a_\Gamma - b_\Gamma$ and $\langle a_\Gamma, n\rangle = 0$ for some vectors $a$ and $b$ have been used.

\begin{figure}[h!]
\begin{center}
  \includegraphics[width=0.9\textwidth]{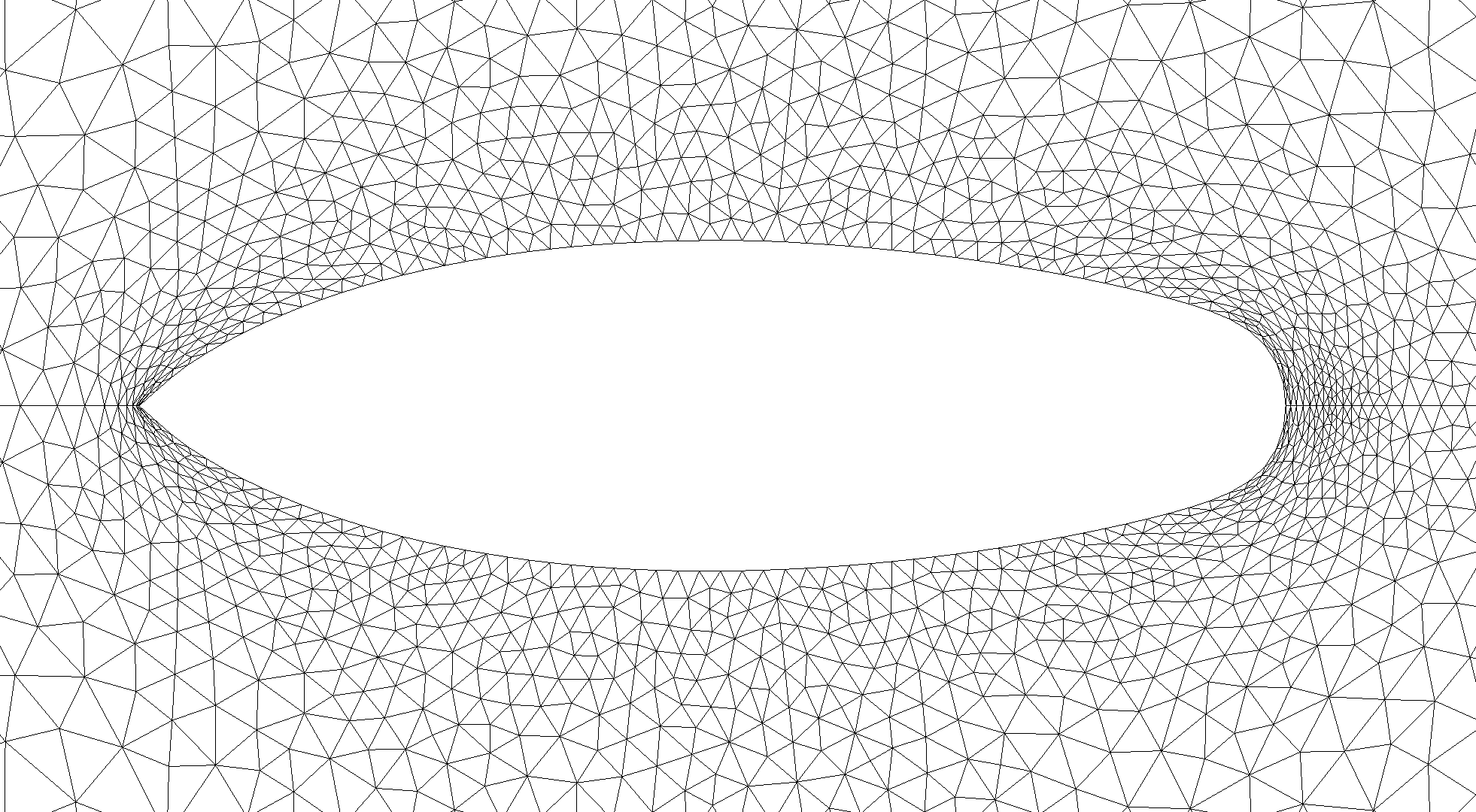}
\end{center}
  \caption{Optimal Shape of the Navier--Stokes problem with surface reparameterization. The volume mesh still has collapsed cells.}\label{fig:SurfaceOkMesh}
\end{figure}

With the same reasoning used when assuming the patches to be independent of the generator in making the transition from~\eqref{eq:MeshBadness} to~\eqref{eq:MeshBadness2}, we also suppose the tangential plane, and thus the normal to it, to be independent of the generator in each repair iteration, which greatly simplifies Equation~\eqref{eq:Repair1} to
\begin{align}\label{eq:Repair2}
x^*_{i, \Gamma} = \frac{\int\limits_{\Gamma_i} \rho(s) s_\Gamma \ \d s}{\int\limits_{\Gamma_i} \rho(s) \ \d s},
\end{align}
which will be used as our tangential mesh repair surface node reparameterization iteration. The resulting improved surface node spacing can be seen in Figure~\ref{fig:SurfaceOkMesh}.
Although the surface mesh point distribution has improved greatly, the issue of collapsed cells still has to be addressed.

\subsection{Eikonal Convection Mesh Deformation}
In principle, the above methodology can also be used to deform the volume mesh mimicking the node spacing aspect of the Delaunay meshing procedure. However, since an even spacing of volume mesh nodes is seldom desired, this approach would first require a reconstruction of the original density $\rho$, which determined the initial node spacing of the original mesh. The most common approach to adapt the volume mesh to a new boundary is to solve a Laplace-like equation, in particular the component wise Laplacian
\begin{alignat*}{3}
 -\epsilon\Delta v &= 0 \quad &\text{ in } &\Omega \\
 v &= g \cdot n \quad &\text{ on } &\partial \Omega,
\end{alignat*}
where $g\cdot n$ is the desired boundary deformation from~\eqref{eq:Update}. This approach, however, often results in collapsed cells and distorted elements if the deformation is too large. Numerous applications within fluid-structure interaction problems have demonstrated that linear elasticity or certain harmonic operators
%\begin{alignat*}{3}
% \div \sigma &= 0 \quad &\text{ in } &\Omega \\
% v &= g \cdot n \quad &\text{ on } &\partial \Omega,
%\end{alignat*}
provide an excellent approach to deform the volume mesh~\cite{MeshDefoAeroElas1, DwightElasDefo, WickMeshMotion1}. Common to those problems is, however, a fairly global and smooth boundary deformation and often those approaches based on elliptic or parabolic PDEs have problems to preserve cell volumes under the highly local deformations occurring in nodal shape optimization. Finally, a very interesting recent approach of volume mesh deformation is to base the PDE to be solved for calculating the nodal displacement vectors on a spectrally equivalent PDE that approximates the pseudo-differential operator nature of the design to state mapping of the optimization problem~\cite{Arian1, schmidt:2562}. That way, extending the deformation to the volume and calculating an approximation to the Newton descent direction happens simultaneously~\cite{PhDSiebenbornGPU}.

As an alternate approach, volume mesh deformation based on a convection-diffusion equation
\begin{alignat}{3}\label{eq:DefoEquation}
 -\epsilon_1\Delta v + \div (v\cdot b) &= 0 \quad &\text{ in } &\Omega \\
 v &= g \cdot n \quad &\text{ on } &\partial \Omega \nonumber
\end{alignat}
is considered here. This idea is to use the convective wind $b := \nabla \epsilon_2$ to project the deformation displacement orthogonal to the surface tangents into the volume mesh. Thus, we follow an inverse ray tracing approach by not calculating an orthogonal projection of every volume node to the design surface, but rather follow an infinite number of orthogonal rays from the design surface into the volume. A convenient way of achieving this is to directly solve the Eikonal equation without rays~\cite{NowackRay1}:
\begin{align}\label{eq:Eikonal}
 h \Delta \epsilon + \Vert \nabla \epsilon \Vert_2^2 = 1 \quad &\text{ in } \Omega,
\end{align}
where $h$ is a small stabilization parameter. When specifying a Dirichlet zero boundary values for $\epsilon$, the numerical value of $\epsilon$ is a non-dimensional travel time to the boundary when locally following the shortest path ray given by $\nabla \epsilon$, which is thus an approximation to the orthogonal projection. Alternatively, the value $\nabla \epsilon$ can be interpreted as a reasonable extension of the surface normal into the volume, which is often required in nodal shape optimization.

\begin{figure}[h!]
\begin{center}
  \includegraphics[width=0.9\textwidth]{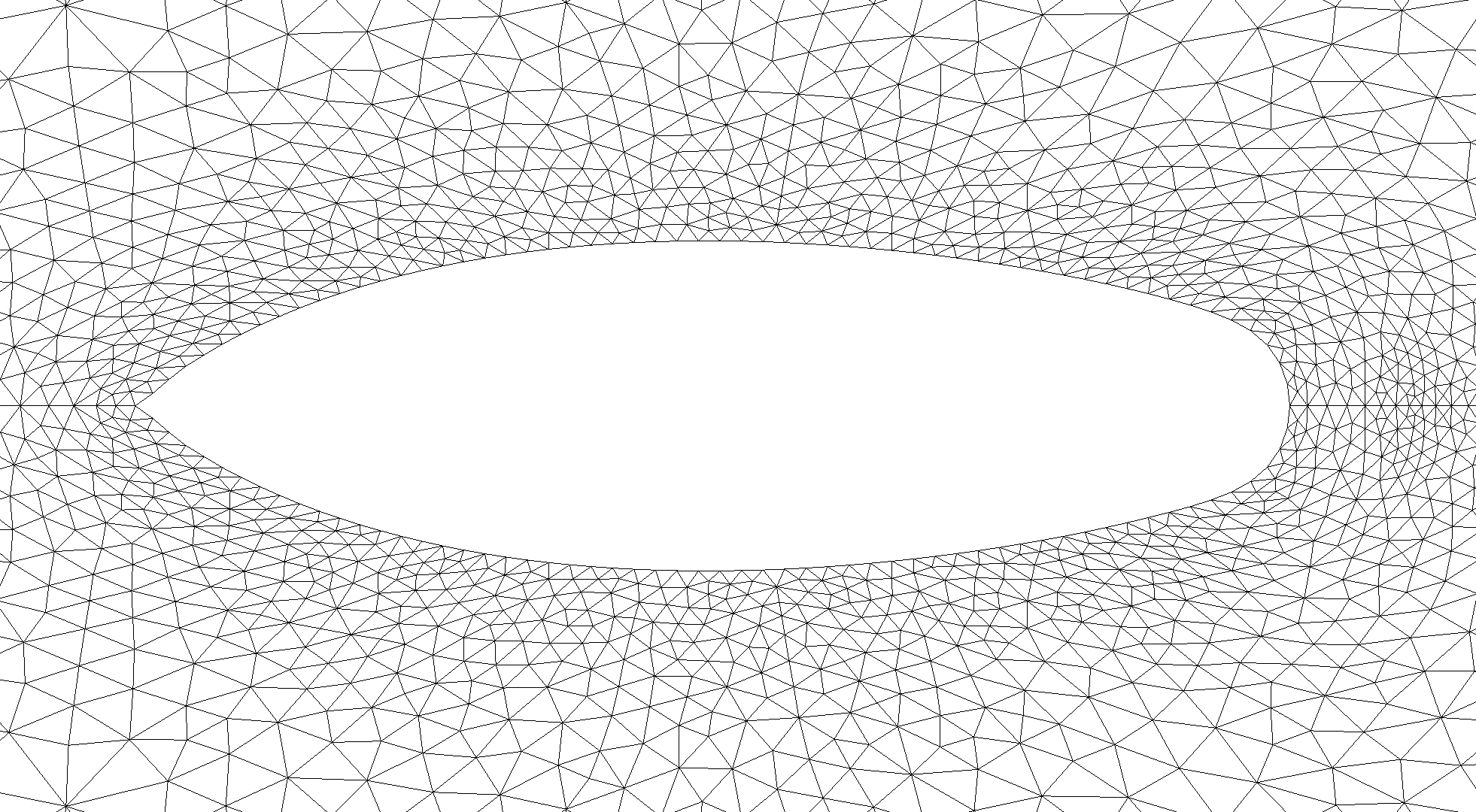}
\end{center}
\caption{The final mesh using both the CVT surface repair and the Eikonal mesh deformation. The surface is nicely parameterized and all volume cells are well shaped. No remeshing is necessary.}\label{fig:Final2D}
\end{figure}

In order to compute the actual mesh deformation in accordance with~\eqref{eq:DefoEquation}, Equation~\eqref{eq:Eikonal} is solved twice, once for $\epsilon_1$ by setting Dirichlet zero boundary conditions on the far-field and once for $\epsilon_2$ by setting Dirichlet zero boundary conditions on the unknown obstacle. Natural boundary conditions are used on the respective other boundaries. Thus, $\epsilon_1(x)$ indicates the travel time, i.e.\ the distance, to the far-field, while $\epsilon_2(x)$ will indicate the distance to the unknown shape to be optimized. The actual displacement vectors are then found by solving
\begin{alignat}{3}
 -\alpha \epsilon_1^2\Delta v + \div (\beta v\cdot \nabla \epsilon_2) &= 0 \quad &\text{ in } &\Omega \\
 v &= g \cdot n \quad &\text{ on } &\Gamma \nonumber \\
 v &= 0 &\text{ on } &\partial \Omega\setminus\Gamma \nonumber
\end{alignat}
where $\alpha$ and $\beta$ are two real parameters to be chosen based on the actual geometry. Diffusivity is scaled down with the square of the distance to the shape changing obstacle, making the convective term dominant in close vicinity. The above equation is solved using 2nd order SUPG-Finite elements. The non-linearity in the Eikonal equation is solved via Newton's method. The resulting mesh is shown in Figure~\ref{fig:Final2D} and is very well behaved.

\section{Application Example: Energy Minimization in a Navier-Stokes Fluid}
\begin{figure}[p]
\begin{center}
\includegraphics[width=0.65\textwidth]{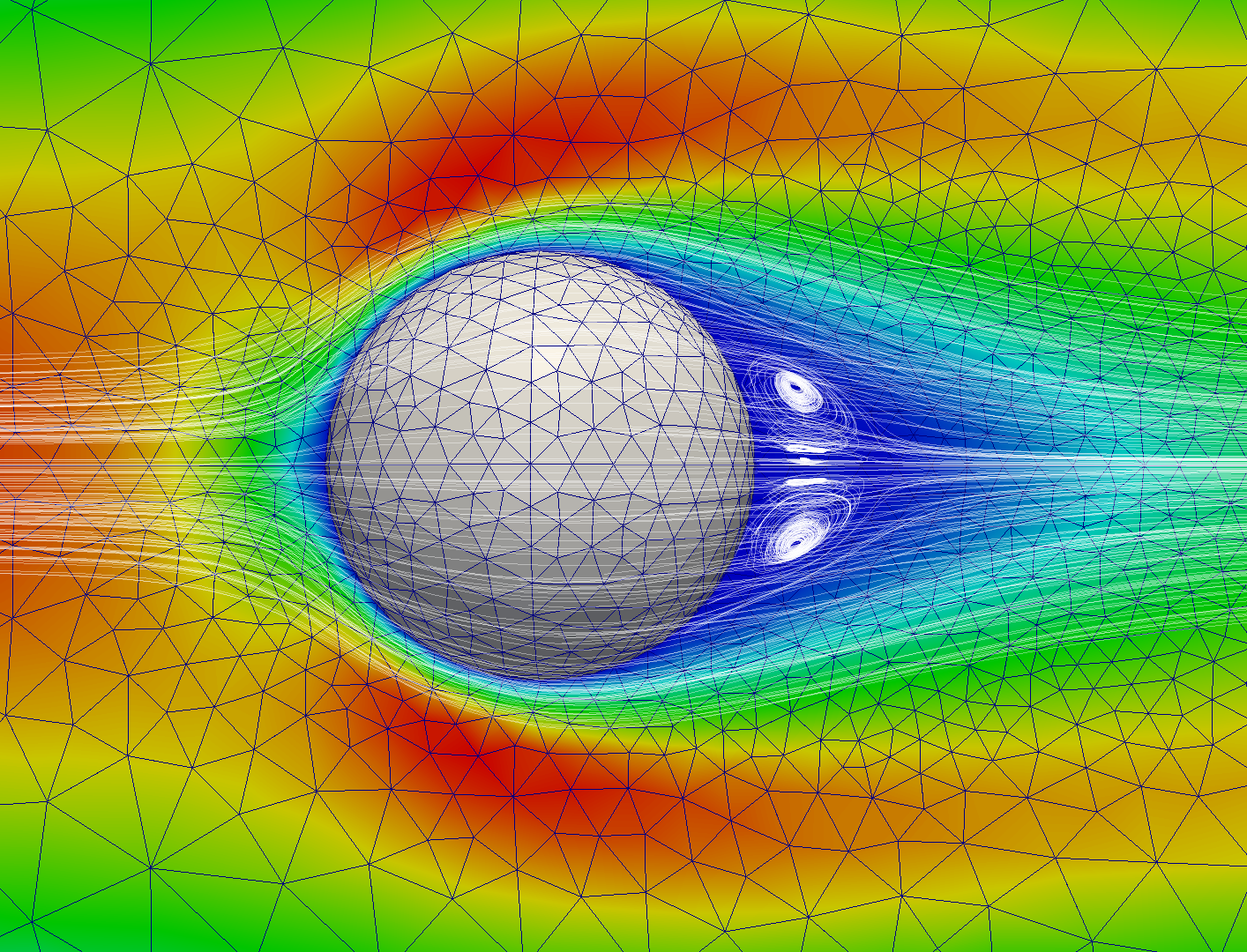} \\
\includegraphics[width=0.65\textwidth]{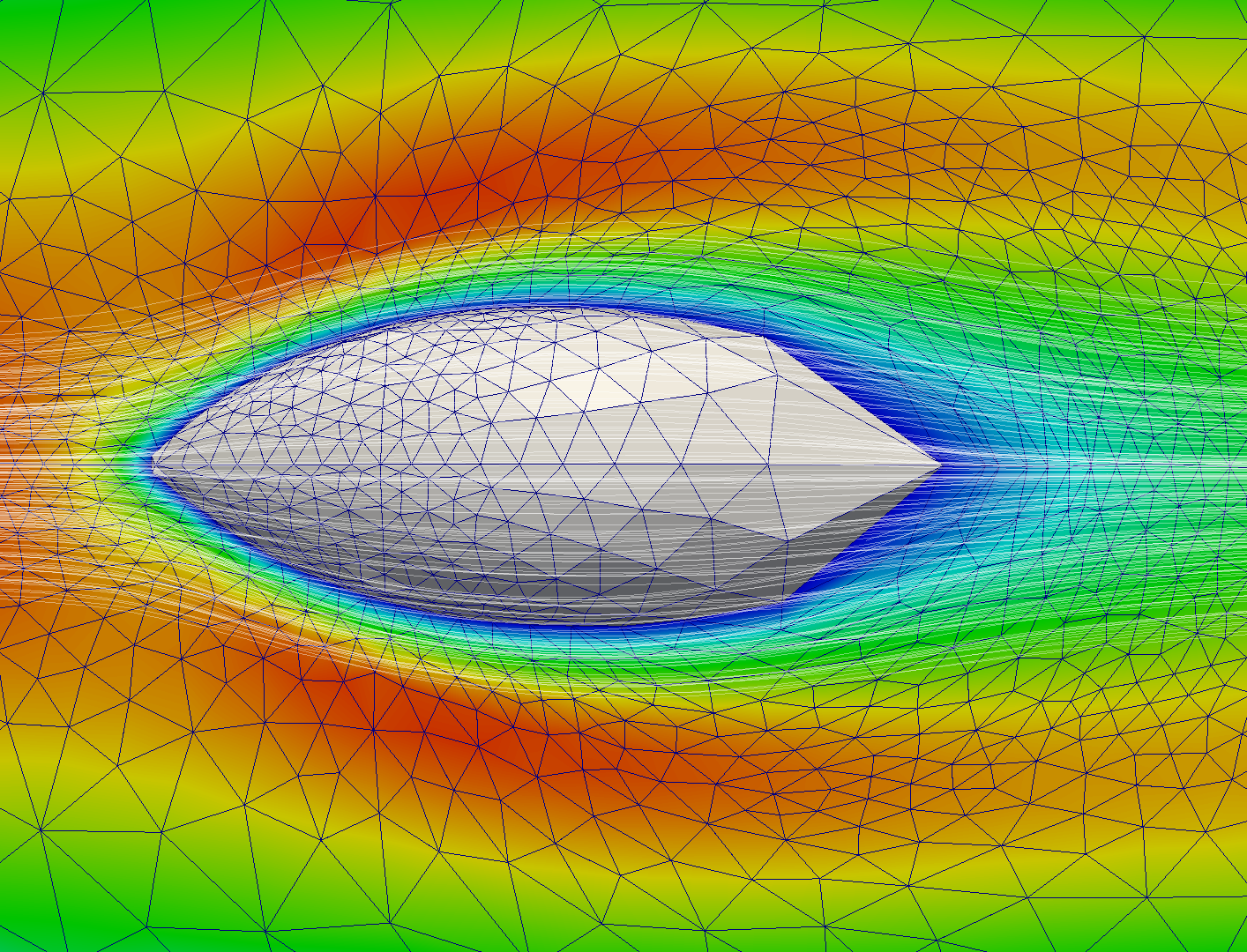} \\
\includegraphics[width=0.65\textwidth]{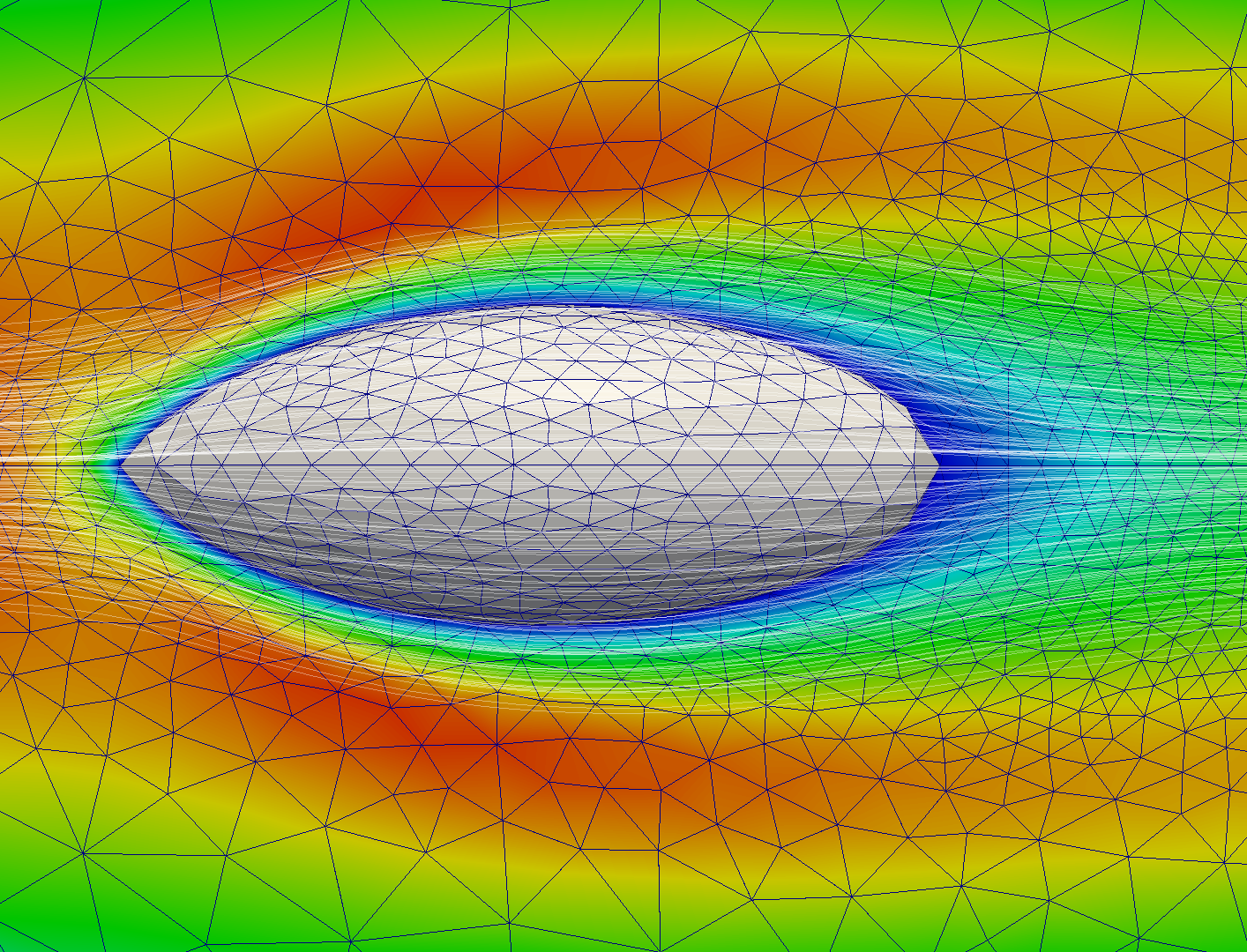}
\end{center}
\caption{Top: Initial 3D geometry. Middle: Naive deformation, no repair, Laplace deformation. Bottom: Full CVT/Eikonal deformation and repair. No remeshing is necessary.}\label{fig:Results3D}
\end{figure}

As an application example, shape optimization in a Navier--Stokes fluid is considered. The problem is well-understood~\cite{BranLinUlbUlb2009, PironneauxX2001, schmidt:2562, SS10GeneralNS} and could be seen as the extension of the DFG simulation benchmark~\cite{DFGBenchmark2D} into shape optimization. Because well-spaced surface node positions are much more difficult so achieve for surfaces $\Gamma$ within a 3D domain $\Omega$, the following problem will be considered in both two and three dimensions.
\allowdisplaybreaks\begin{align}
  &\begin{aligned}\label{eq:Objective}
    \min\limits_{(\velo, \p, \Omega)} J(\velo, \p, \Omega) &:= \int\limits_\Omega \dvis \sum\limits_{i,j=1}^3 \left(\frac{\partial \velo_i}{\partial x_j}\right)^2 \ dA
    \end{aligned} \\
 &\mbox{subject to}\nonumber \\
 &\begin{aligned}
    -\dvis \Delta \velo + \dens \velo \nabla \velo + \nabla \p &= 0 &&\text{in}\quad \Omega \\
    \div \velo & = 0 &&\\
    \velo &= \velo_+ &&\text{on}\quad \Gamma_+\\
    \velo &= 0 &&\text{on}\quad \Gamma_0 \\
    \p \n - \dvis \frac{\partial \velo}{\partial \n} &= 0 &&\text{on}\quad \Gamma_- \\
    \operatorname{vol} &= V_0 \\
    \operatorname{geo}(\Omega) &= 0,
  \end{aligned}
\end{align}
where $\operatorname{vol}$ denotes the volume and $\operatorname{geo}$ summarizes additional geometric constraints such as fixing the centroid of $\Omega$ and preserving the symmetry planes. Especially the latter is quite important, as otherwise an incorrect discrete normal at sharp edges~\cite{Lozano2012} and general finite element approximation errors~\cite{MorinAdaptive1} tend to result in a problematic behavior at spurious or physical geometric singularities. Not considering the volume and geometric constraints, the shape derivative of the above problem is given by
\begin{align}
  &\begin{aligned}
    dJ(\velo, \p, \Omega) =\int\limits_{\Gamma_0} \langle\V,\n\rangle\left[-\dvis\sum_{i=1}^3\frac{\partial \lambda_i}{\partial \n}\frac{\partial \velo_i}{\partial \n} + \left(\frac{\partial \velo_i}{\partial \n}\right)^2\right]\ dS,
    \end{aligned}\label{eq:NavierStokesShapeGradient}
 \end{align}
where $(\lambda, \lambda_p)$ solves the adjoint Navier--Stokes equations
\begin{align}
 &\begin{aligned}\nonumber
    -\dvis \Delta \lambda_i - \dens \sum_{j=1}^3 \left(\frac{\partial \lambda_j}{\partial x_i} \velo_j + \frac{\partial \lambda_i}{\partial x_j}\velo_j\right) - \frac{\partial \lambda_\p}{\partial x_i} &= -2\dvis\Delta \velo_i && \text{in} \quad \Omega\\
    \div \lambda &= 0& \\
    \lambda &= 0 && \text{on} \quad \Gamma_+\\
    \lambda &= 0&&\text{on} \quad \Gamma_0\\
    \dvis \frac{\partial \lambda_i}{\partial \n} + \dens\left(\sum_{j=1}^3 \lambda_j\velo_j\n_i + \lambda_i\velo_j\n_j\right) + \lambda_\p\n_i&= 0 &&\text{on} \quad \Gamma_-.
  \end{aligned} 
\end{align}
Both primal and dual Naiver--Stokes equations are implemented in Python FEniCS. The non-linearity of the primal problem is either solved using a Newton-Iteration or by successively solving Oseen problems stemming from a Picard-linearization. From a linear-algebra perspective, the resulting saddle-point problem is either directly factorized using the parallel sparse direct solver ``mumps'' or solved iteratively using a continuous pressure Schur-complement iteration inside an implicit Euler time stepping to steady state~\cite{TurekBook}. Initial and optimal shape with and without the Eikonal mesh deformation and repair step are shown in Figure~\ref{fig:Results3D}.

In 2D, the domain is a rectangle with a half-axis of $1.5$ around a half-circle of radius $0.5$ at the coordinate origin. The three dimensional geometry is created by  rotating this setup around the $x$-axis by $90$ degrees. This half-domain, or quarter-domain in 3D, is discretized using an unstructured mesh. The actual computational mesh is afterwards created by mirroring these half-meshes, thereby creating a perfectly symmetric mesh that can still be used to represent non-symmetric solutions, e.g. separated flows and vortex streets. The primal and adjoint Navier--Stokes equation is discretized using LBB-stable P2P1-Taylor Hood elements, leading to $131,790$ unknowns for the primal and dual velocities each and $5,728$ unknowns for each primal and dual pressure in 3D. For the two dimensional test-cases we arrive at $17,324$ and $2,221$ unknowns respectively. The Eikonal equation and the convection-diffusion equation appearing during the mesh adaptation step are also solved using 2nd order SUPG-stabilized finite elements. All PDEs are discretized using the Python implementation of FEniCS, which also automatically generates the Newton solver for the Eikonal equation.

The reduced Hessian of the problem is approximated by a Laplace-Beltrami operator, such that a smoothed gradient descent scheme is employed. For more details see~\cite{schmidt:2562}. Volume and miscellaneous geometric constraints are enforced using a simple projection during each shape update. The actual algorithm is described in Algorithm~\ref{algo:1} in all detail.
\begin{algorithm}[h!]
\begin{algorithmic}
\STATE $q_k \gets$ solution of state equation in current domain $\Omega_k$
\STATE $\lambda_k \gets$ solution of adjoint equation in current domain $\Omega_k$
\STATE Calculate boundary displacement field $d_1\gets g(q_k,\lambda_k)$ as in Equation~\eqref{eq:Update}
\STATE Calculate approximate Newton update: $d_2 \gets$ solution of $(\delta \Delta_\Gamma + I) d_2 = -d_1$
\STATE Project $d_2$ to geometric constraints $\operatorname{geo}$ and volume constraint $\operatorname{vol}$
\STATE Calculate intermediate next boundary $\partial \Omega_{k+\frac{1}{2}} = \{x + d_2(x)\cdot n(x) : x \in \partial \Omega_k\}$
\LOOP
 \STATE Calculate intermediate surface node positions $x_{i,\Gamma}^*$ by solving \eqref{eq:Repair2}.
 \STATE Calculate tangential offset and badness $\tau(x) \gets x_{i,\Gamma}^* - (x + d_2(x)\cdot n(x))$
 \IF {$\Vert \tau \Vert \leq \epsilon$}
   \STATE break
 \ENDIF
\STATE $\partial \Omega_{k+\frac{1}{2}} \gets \{x + d_2(x)\cdot n(x) + \tau : x \in \partial \Omega_k\}$
\ENDLOOP
\STATE Define final normal and tangential boundary displacement $d_3 := d_2 \cdot n + \tau$
\STATE $\epsilon_1 \gets$ solution of Eikonal~\eqref{eq:Eikonal} in $\Omega_k$ with Dirichlet zero on far-field
\STATE $\epsilon_2 \gets$ solution of Eikonal~\eqref{eq:Eikonal} in $\Omega_k$ with Dirichlet zero on design
\STATE Define $b := \nabla \epsilon_2$ ray tracing directions in $\Omega_k$
\STATE Calculate full deformation in $\Omega$, i.e.\ $d_4 \gets $ solution of~\eqref{eq:DefoEquation}
\STATE Update domain $\Omega_{k+1} = \{x + d_4(x) : x \in \Omega_{k}\}$
\IF {$\Omega_{k+1}$ fulfills convergence}
	\STATE exit
\ELSE
	\STATE loop
\ENDIF
\STATE
\end{algorithmic}
\caption{Nodal shape optimization scheme with CVT step and Eikonal deformation.}\label{algo:1}
\end{algorithm}

It is worth noting that, contrary to the aerodynamic force minimization using boundary integrals, the energy objective, Equation~\eqref{eq:Objective}, scales with the domain size $\Omega$, which makes an absolute comparison of objective difficult. Nevertheless, for the results presented here, the objective is reduced from $266.1282$ to $244.9740$ in the 3D case, and from $195.6071$ to $131.3547$ in the 2D situation, a relative decrease comparable to~\cite{schmidt:2562}. The initial and respective final shapes are shown in Figure~\ref{fig:Results3D}. With the CVT tangential movement and Eikonal mesh deformation, the optimal shape is nicely parameterized and the whole optimization works successfully, maintaining well-behaved volume cells without re-meshing.

\section{Conclusions}
A two level mesh repair and mesh deformation strategy for large nodal shape optimization is considered and exemplified using the energy minimization in an incompressible Navier--Stokes fluid. A well-parameterized surface is maintained by introducing a mesh quality constraint into the nodal boundary update based on minimizing the tangential projection of the mesh badness criterion stemming from the centroidal Voronoi tessellation, the mesh dual of the Delaunay triangulation. The methodology is useable as a surface mesh deformation scheme, rather than a re-meshing approach, because the patch definition is not changed within each quality preservation step, meaning the cell connectivity is not changed. Furthermore, the tangential repair step is in the kernel of the shape derivative, thus a reparameterization has no impact on optimality in a continuous interpretation.

In a second step, the boundary deformation is propagated into the volume mesh by solving a convection--diffusion equation with wind direction based on a ray tracing approach, which is achieved by solving an Eikonal equation. The methodology is exemplified by minimizing the energy dissipation in a Navier--Stokes fluid completely without remeshing in both two and three dimensions.

\bibliography{MyBib}{}
\bibliographystyle{plain}

\end{document}